\documentclass{amsart}

\newtheorem{thm}{Theorem}[section]

\theoremstyle{definition}

\theoremstyle{remark}

\newtheorem*{rem}{Remark}

\numberwithin{equation}{section}

\newcommand{\slt}{\mathfrak{sl}(3)}
\newcommand{\sln}{\mathfrak{sl}(n)}
\newcommand{\g}{\mathfrak{g}}
\newcommand{\frakh}{\mathfrak{h}}

\newcommand{\ts}{\otimes}
\newcommand{\R}{\mathfrak R}
\newcommand{\Or}{{\mathcal O}_{r}}
\newcommand{\OR}{{\mathcal O}_{R}}
\newcommand{\Og}{{\mathcal O}_{\gamma}}

\DeclareMathOperator{\im}{Im}
\DeclareMathOperator{\ad}{ad}
\begin{document}


\title[Boundary $R$-matrices]
{Boundary solutions of the quantum Yang-Baxter equation and
solutions in
three dimensions}
\author{Murray Gerstenhaber}
\address{Department of Mathematics, University of Pennsylvania,
Philadelphia, PA 19104-6395}
\email{mgersten@sas.upenn.edu}

\author{Anthony Giaquinto}
\address{Department of Mathematics, Texas A \& M University,
College Station, TX 77843-3368}
\email{tonyg@math.tamu.edu}
\thanks{The second author was supported in part by a grant from
the National
Security Agency}

\date{October 21, 1997}



\begin{abstract}
Boundary solutions to the quantum Yang--Baxter (qYB) equation are
defined to be those in the boundary of (but not in) the variety
of solutions to the ``modified'' qYB equation, the latter being
analogous to the modified classical Yang--Baxter (cYB) equation.
We construct, for a large class of solutions $r$ to the modified
cYB equation, explicit ``boundary quantizations'', i.e., boundary
solutions to the qYB equation of the form $I+tr+ t^2r_{2} + \dots$. In
the last section we list and give quantizations for all classical
$r$-matrices in $\slt \wedge \slt$.
\end{abstract}

\maketitle


\section{Introduction}
Solutions to the quantum Yang-Baxter (qYB) equation,
so-called ``$R$-matrices,'' have many applications in mathematical
physics and are essential to the theory of quantum groups, but
explicit solutions are scarce. Finding directly all constant
$n^{2} \times n^{2}$ matrix solutions reduces to solving $n^{6}$
homogeneous cubic equations in $n^{4}$ variables, a task feasible
only in the case $n=2$, see \cite{H}. Usually the search is
restricted to one for those $R$-matrices which are quantizations
(in the sense below) of ``$r$-matrices'', i.e., solutions to the
classical Yang-Baxter (cYB) equation, but even these $R$ are
somewhat elusive. The problem is two-fold -- first find all $r$-matrices, 
then quantize them. In \cite{GG}, we introduced the
notion of a boundary solution to the classical Yang-Baxter
equation in the hope that there may be some reasonable
classification for these. Section 2 introduces the natural
extension of this concept to that of a boundary solution to the
quantum Yang-Baxter equation. For a large class of boundary
classical $r$-matrices we explicitly construct in Section 3
associated boundary quantizations. Finally, Section 4 contains a
complete classification (up to conjugacy) of all $R$-matrices
which are quantizations of classical $r$-matrices in
$\slt \wedge \slt$. While not all of these $R$ are constructed as
boundary solutions, in a sense the most interesting ones are. It
remains an open question to determine precisely which ones
actually are boundary $R$-matrices.

 \section{The quantum Yang-Baxter equations}
     Let $V$ be an $n$-dimensional vector space over a field $k$.
A quantum $R$-matrix is an operator $R:V\ts V\rightarrow V\ts V$
satisfying the quantum Yang-Baxter equation
  $$R_{12}R_{13}R_{23}= R_{23}R_{13}R_{12}$$
where as usual $R_{ij}$ denotes $R$ operating in the $i$ and $j$
factors of the triple tensor product $ V^{\otimes3} = V\otimes
V\otimes V$. We define the {\it modified} qYB equation (MQYBE) to
be
$$R_{12}R_{13}R_{23}- R_{23}R_{13}R_{12} = \lambda ((123)R_{12}-
(213)R_{23})$$
where $\lambda$ is a non-zero element of the coefficient field
and $(123)$, $(213)$ denote the corresponding permutations of the tensor
factors of $V^{\otimes3}$. Correspondingly, we call $R$ a
modified $R$-matrix if it satisfies the MQYBE. A solution to
either equation is {\it unitary} if $R_{21}= R^{-1}$.
 \begin{rem} Modified $R$-matrices are closely connected to
Hecke symmetries, for if $R$ is a modified unitary $R$-matrix
then  $((12)R-\lambda \, I)/\sqrt {(1-\lambda ^2)}$ satisfies the
braid and Hecke relations with $q=(1-\lambda)/\sqrt {1-\lambda
^2}$, see \cite{GGS} for a detailed
discussion.
 \end{rem}

  A natural problem is to determine the structure of the spaces
of solutions to both the original and modified quantum
Yang-Baxter equations. It is clear that if $R$ satisfies either
then so does any scalar multiple of $R$, so it is meaningful to
view  $R\in {\mathbb P}(M_{n^{2}}(k))$. While little is known in
general, the following simple result is helpful in determining
many solutions.
 \begin{thm} \label{thm1}
 Let $\mathfrak R$, respectively $\mathfrak R'$, denote the
subsets of
 ${\mathbb P}(M_{n^{2}}(k))$ consisting of solutions to the
 quantum Yang-Baxter equation, respectively
 modified quantum Yang-Baxter equation,
 and let $\overline {\mathfrak R'}$ be the Zariski closure of
 ${\mathfrak R'}$. Then
 \begin{enumerate}
 \item $\mathfrak R$ is a closed variety.
 \item $\overline {\mathfrak R'} - \mathfrak R' \subset \mathfrak
R.$
 \item $\mathfrak R'$ is quasi-projective, i.e. it is an open
subset of its
 closure.
 \end{enumerate}
 \end{thm}
 \begin{proof} \mbox { }\\
 \begin{enumerate}
 \item Each entry of $R_{12}R_{13}R_{23}- R_{23}R_{13}R_{12}$ is
a homogeneous
 cubic polynomial in the entries of $R$ and $\R$, being the locus
of common zeroes of these polynomials, is therefore a closed
variety.
 \item Define maps $\phi :({\mathbb P}(M_{n^{2}}(k)) -
\R)\rightarrow
 {\mathbb P}(M_{n^{3}}(k))$ with $\phi (R)= R_{12}R_{13}R_{23}-
R_{23}R_{13}R_{12}$
 and $\psi : ({\mathbb P}(M_{n^{2}}(k)) - I)\rightarrow
 {\mathbb P}(M_{n^{3}}(k))$ with $\psi (R)=((123)R_{12}-
(213)R_{23})$.
 Since $\R '=\phi ^{-1}(\im (\phi) \cap \im (\psi))$ and
 $\im (\phi) \cap \im (\psi)$ is closed, it follows that $\R '$
is closed in
 ${\mathbb P}(M_{n^{2}}(k)) - \R$. Thus $\overline {\mathfrak R'}
- \mathfrak R'
 \subset \mathfrak R.$
 \item To show that $\R '$ is open in $\overline {\R '}$, it
suffices to verify
 that $\overline {\R '} -\R '$ is closed in $\overline {\R '}$.
Since
 $\overline {\R '} -\R ' =  \overline {\R '} \cap \R $ and
 both $\overline {\R '}$ and $\R $ are closed, it follows that
 $\overline {\R '} -\R '$ must also be closed.
 \end{enumerate} \end{proof}
 In light of Theorem \ref{thm1}, we call an $R\in \overline {\R
'} -\R '$ a {\it boundary solution} to the QYBE or boundary $R$-
matrix.
\section{Construction of boundary solutions}
In \cite{GG} we introduced boundary solutions of the classical
Yang-Baxter equation and gave some procedures for constructing
them. In this section we ``quantize'' some of those.

Suppose that $R$ is a modified unitary quantum $R$-matrix of
the form $R=I +tr +O(t^{2})$. Then the ``(semi)classical limit''
$r$ is a skew modified classical $r$-matrix, that is,
$[r_{12},r_{13}] + [r_{12},r_{23}] + [r_{13},r_{23}]$ is a
non-zero $\sln$-invariant. Similarly, the classical limit of a
quantum unitary $R$-matrix is a skew classical $r$-matrix. In
either case we call $R$ a ``quantization'' of $r$.
Now let $r \in \sln \wedge \sln$ be a modified $r$-matrix, and
suppose that $R$ is a modified unitary $R$-matrix
with classical limit $r$. Let $\Or$ and $\OR$ denote,
respectively, the
$SL(n)$-orbits of $r$ and $R$ and let
$\overline {{\mathcal O}}_{r}$ and
$\overline {{\mathcal O}}_{R}$ be their closures.
For any $\tilde{r}\in \overline {{\mathcal O}}_{r}
-\Or$,
we seek an $\widetilde{R}\in \overline {{\mathcal O}}_{R}- \OR$
which quantizes $\tilde{r}$.
In \cite{GG} it was shown that such an $\tilde{r}$ is a boundary
classical $r$-matrix so what we are really looking for is an
associated boundary $R$-matrix. There is presently no complete
answer to this question, although a conjectured explicit form for
the quantization of all modified classical $r$ does exist, see
\cite{GGS}. (It is known that every such $r$ does admit
a quantization.) We can however give a positive answer in many
cases.
\begin{rem} Given an $r$, the quantizations we produce will be of
the form
$R=1+2tr +O(t^{2})$ which is, course, equivalent to finding an
$R$ with
classical limit $r$. This will be particularly useful in the next
section
as we can eliminate powers of $1/2$ in the $R$-matrices.
\end{rem}

Set $\gamma = \sum _{i<j} e_{ij} \wedge e_{ji}$, this is
commonly called the Drinfel'd-Jimbo modified classical $r$-matrix
since it is the classical limit of the standard quantization of
$SL(n)$. It was shown in \cite{GGS} that the appropriate
quantization of $\gamma$ is $\exp(2t\gamma )$. This is a modified
unitary quantum $R$-matrix and the corresponding Hecke symmetry
coincides with
that of the classic $R$-matrix associated with ${\mathcal
O}_{q}(SL(n))$.
\begin{thm} \label{thm2}
Let $\gamma '\in \overline {{\mathcal O}}_{\gamma} -\Og$.
Then $\exp(2t\gamma ')$ is a
boundary solution of the quantum Yang-Baxter equation.
\end{thm}
\begin{proof}
Since $\gamma '\in \overline {{\mathcal O}}_{\gamma}$, it follows
that
$\exp(2t\gamma ')\in \overline {\mathcal {O}}_{\exp(2t\gamma )}
\subset \overline {\R}'.$
Now because $\gamma '$ is a classical $r$-matrix,
$\exp(2t\gamma ')\notin \R'$ and so by Theorem \ref{thm1} it
follows that
$\exp(2t\gamma ')$ is a boundary solution to the quantum
Yang--Baxter equation.
\end{proof}

As is well known, $\gamma + \beta $ is a modified $r$-matrix
for any $\beta \in \frakh \wedge \frakh$.
These are the classical limits of the standard multi-parameter
$R$-matrices. The modified quantum $R$-matrix which quantizes
$\gamma + \beta $ is $\exp(t\beta )\cdot \exp(2t\gamma)\cdot
\exp(t\beta)$, see \cite{GGS}. For a large class of boundary
classical $r$-matrices in
$\overline {\mathcal {O}}_{\gamma + \beta}-\mathcal {O}_{\gamma +
\beta}$
the same type quantization is possible.
\begin{thm} \label{thm3}
Suppose $x\in \sln$ and suppose
$$\exp(v\ad x)(\gamma +\beta) =
(\gamma + \beta) + v(\gamma_{1}+\beta_{1})+ \cdots +
v^{d}(\gamma_{d}+\beta_{d}).$$
Then $\gamma_{d}+\beta_{d}\in
\overline {\mathcal {O}}_{\gamma + \beta}-\mathcal {O}_{\gamma +
\beta}$
and $\exp(t\beta_{d} )\cdot \exp(2t\gamma_{d})\cdot
\exp(t\beta_{d})$ is
a boundary $R$-matrix which quantizes $\gamma_{d}+\beta_{d}$.
\end{thm}

\begin{proof}
According to \cite{GG},
$\gamma_{d} +\beta_{d}\in
\overline {\mathcal {O}}_{\gamma + \beta}-\mathcal {O}_{\gamma +
\beta}$ and
so is a boundary classical $r$-matrix. Set
$M=\exp(t\beta )\cdot \exp(2t\gamma)
\cdot\exp(t\beta)$.
Then $\exp(t\beta_{d} )\cdot \exp(2t\gamma_{d})\cdot
\exp(t\beta_{d})\in
\overline {\mathcal {O}}_{M}\subset \overline {\R}'$.
But since its classical limit
is an $r$-matrix,
$\exp(t\beta_{d} )\cdot \exp(2t\gamma_{d})\cdot
\exp(t\beta_{d})\notin
\R'$ and so by Theorem \ref{thm1} it must be a boundary solution
to the
quantum Yang-Baxter equation.
\end{proof}

\section{Quantization of the classical $r$-matrices in
$\slt \wedge \slt$}
We now consider the quantizations of all classical
$r$-matrices in $\slt \wedge \slt$.
The quantizations of the modified classical $r$-matrices in $\slt
\wedge \slt$
are well known. Up to automorphism of
$\slt$, there are two distinct types of such $r$,
the one-parameter standard family
$$r_{\lambda}=e_{12}\wedge e_{21}+ e_{13}\wedge e_{31}+e_{23}
\wedge e_{32} + \lambda (e_{11}-e_{22})\wedge (e_{22}- e_{33})$$
and what we have called the Cremmer-Gervais solution
$$r_{CG}= r_{1/3}+ 2e_{12}\wedge e_{32}.$$
The modified $R$-matrix associated with $r_\lambda$ was discussed
in the
previous section and
the Cremmer-Gervais modified $R$ is
$$q^{-\beta} \left\{\exp (2t\gamma)+2\sin
(t)(q^{1/2}e_{12}\otimes e_{32}
-q^{-1/2}e_{32}\otimes
e_{12})\right\}q^{-\beta}$$
with $\beta =(1/3) (e_{11}-e_{22})\wedge (e_{22}- e_{33})$ and
 $q=\sec(t) -\tan (t)$, see \cite{GGS}.

Now let us focus on the classical $r$-matrices in $\slt \wedge
\slt$.
Using the homological interpretation (due to Drinfel'd) of the
classical Yang-Baxter equation, Stolin \cite{S} has given a
description of classical $r$-matrices in terms of ``quasi-Frobenius'' 
Lie algebras. A Lie algebra $\mathfrak f$ is
quasi-Frobenius if there is a non-degenerate map $\phi
:{\mathfrak f}\wedge {\mathfrak f}\rightarrow k$. The inverse
matrix of $\phi$ is then a classical $r$-matrix. Now if
$\mathfrak g$ is a simple Lie algebra, and $r\in {\mathfrak
g}\wedge \mathfrak g$ is a classical $r$-matrix, then
there is a unique proper subalgebra $\mathfrak f$ of $\g$ for
which $r\in {\mathfrak f}\wedge {\mathfrak f}$ and is non-degenerate.
 Call $\mathfrak f$ the {\it carrier} of $r$.
Exploiting the relatively small dimension, Stolin gave a complete
list of quasi-Frobenius subalgebras of $\slt$ and computed enough
of their cohomology to determine the number of classical $r$-matrices
 each carries up to automorphism. However, no explicit
form of these $r$ nor of their quantizations was given in
\cite{S}.

The rest of this note is devoted to listing explicitly all the
classical $r$-matrices in $\slt \wedge \slt$ and finding unitary
quantizations of each. If ${\mathfrak f}$ is quasi-Frobenius then
$r({\mathfrak f})$ will denote a classical $r$-matrix carried by
${\mathfrak f}$ and $R({\mathfrak f})$ will be an $R$-matrix
which quantizes $r({\mathfrak f})$.

Not all $r$-matrices in $\slt \wedge \slt$ seem to be boundary
solutions to the classical Yang-Baxter equation so other
techniques besides those in
Theorems \ref{thm2} and \ref{thm3} are needed to quantize them.
The following theorem from \cite{GZ} is useful here.
\begin{thm} \label{thm4}
Let ${\mathfrak s}$ be the two-dimensional Lie algebra generated
by
$H$ and $E$ with relation $[H,E]=E$. Set $H^{\langle d \rangle}=
H(H+1)(H+2)+\cdots (H+(d-1))$. Let
$F_{m}=\sum _{i=0}^{m}=(-1)^{i}\binom{m}{i}E^{m-i}H^{\langle i
\rangle}\otimes
E^{i}H^{\langle m-i \rangle}$ and let $F=\sum
_{i=0}^{\infty}(t^{i}/i!)F_{m}$.
Then $F_{21}^{-1}F$ is a universal quantization of $2(E\wedge
H)$.
\end{thm}
Specializing the universal $R$ to any finite dimensional
representation of
${\mathfrak s}$ gives an $R$-matrix.
Another technique which sometimes
 quantizes non-boundary $r$-matrices is to simply exponentiate
them. We are not sure which class of $r$-matrices has the
property that
the exponential map provides an $R$-matrix, but, fortunately, for
those
$r$ in $\slt \wedge \slt$ which are not covered by Theorems
\ref{thm2},
\ref{thm3},
or \ref{thm4} it turns out that the exponential map works.

Up to automorphism, there are 10 classes of quasi-Frobenius
subalgebras of $\slt$. Their dimensions must be even.
We will for the most part  follow the numbering and notation of
these Lie algebras used in \cite{CP}, Section 3.1.

\begin{enumerate}
\item[(i)]
The only six dimensional subalgebra of $\slt$ is
$${\mathfrak p}=\left(
\begin{matrix}*&*&*\\0&*&*\\0&*&*\end{matrix}\right)$$
and $r({\mathfrak p})=(2e_{11}-e_{22}-e_{33})\wedge e_{12}+ (
e_{11}+
 e_{22}-2e_{33})\wedge e_{23}
+3e_{13}\wedge e_{32}$ is the unique classical $r$ with carrier
${\mathfrak p}$,
see \cite{GG}.
Let $Q$ be the Cremmer-Gervais modified $R$. Let $x=-e_{12}-
(1/2)e_{23}$ and set
$R({\mathfrak p})=\lim_{t\to 0}\exp\left(\frac{3m}{t}\ad
(x)\right)\cdot Q.$
Then
$$R({\mathfrak p})=\left(
\begin{matrix}
1 & 2\,m & m^{2} &  - 2\,m & 2\,m^{2} & m^{3} & m^{2} &  - m^{3}
 &  - 2\,m^{4} \\
0 & 1 & m & 0 & m & m^{2} & 0 & m^{2} &  - 2\,m^{3} \\
0 & 0 & 1 & 0 & 0 & m & 0 & 3\,m &  - 2\,m^{2} \\
0 & 0 & 0 & 1 &  - m & m^{2} &  - m & m^{2} & 2\,m^{3} \\
0 & 0 & 0 & 0 & 1 & m & 0 &  - m & 2\,m^{2} \\
0 & 0 & 0 & 0 & 0 & 1 & 0 & 0 & 2\,m \\
0 & 0 & 0 & 0 & 0 &  - 3\,m & 1 &  - m &  - 2\,m^{2} \\
0 & 0 & 0 & 0 & 0 & 0 & 0 & 1 &  - 2\,m \\
0 & 0 & 0 & 0 & 0 & 0 & 0 & 0 & 1
\end{matrix}\right) $$
is a boundary $R$-matrix which quantizes $r({\mathfrak p})$.

\item[(ii)]
The four dimensional subalgebra
$${\mathfrak r}=\left(
\begin{matrix}*&*&*\\0&*&0\\0&0&*\end{matrix}\right)$$
carries a one-parameter family of classical $r$-matrices.
If $x=(3a/2)e_{12}+(3/2)e_{13}$ then
$\exp(t\ad x)\cdot r_{1/3} = r_{1/3}+t[x,r_{1/3}]$ and so
 $[x,r_{1/3}]= d_{1}\wedge e_{12} + d_{2}\wedge e_{13}$ where
 $d_{1}= a(-e_{11}+2e_{22}-e_{33})+3e_{23}$ and $d_{2} =
-2e_{11}+e_{22}+e_{33}$ is a boundary $r$-matrix.
 Now it is easy to check that the Lie algebra spanned
by $d_{1}, d_{2}, e_{12}$ and $e_{13}$ is conjugate to $\mathfrak
r$ and so it
makes sense to denote
$[x,r_{1/3}]$ as $r({\mathfrak r})$. Now applying Theorem
\ref{thm3}, we get
the boundary matrix
$$R(\mathfrak r)=\left(
\begin{matrix}
1 &  - a\,t &  - 2\,t & a\,t & 2\,a^{2}\,t^{2} & a\,t^{2} & 2\,t
 & a\,t^{2} & 2\,t^{2} \\
0 & 1 & 0 & 0 &  - 2\,a\,t &  - 3\,t & 0 &  - t & 0 \\
0 & 0 & 1 & 0 & 0 & a\,t & 0 & 0 &  - t \\
0 & 0 & 0 & 1 & 2\,a\,t & t & 0 & 3\,t & 0 \\
0 & 0 & 0 & 0 & 1 & 0 & 0 & 0 & 0 \\
0 & 0 & 0 & 0 & 0 & 1 & 0 & 0 & 0 \\
0 & 0 & 0 & 0 & 0 & 0 & 1 &  - a\,t & t \\
0 & 0 & 0 & 0 & 0 & 0 & 0 & 1 & 0 \\
0 & 0 & 0 & 0 & 0 & 0 & 0 & 0 & 1
\end{matrix}\right) $$
which is clearly a quantization of $r({\mathfrak r})$.

\item[(iii)]
The two-parameter family of four-dimensional subalgebras of
$\slt$
$${\mathfrak q}_{a,b,c}=\left(
\begin{matrix}a&0&0\\0&b&0\\0&0&c\end{matrix}\right)+
\left(
\begin{matrix}0&*&*\\0&0&*\\0&0&0\end{matrix}\right)$$
also carries classical $r$-matrices.
For $(a,b,c)=(0,1,-1)$ or $(1,1,-2)$, there two classical $r$-matrices
carried by ${\mathfrak q}_{a,b,c}$, for
$a=c$ there are no classical $r$-matrices carried by
${\mathfrak q}_{a,b,c}$, and the remaining types each carry a
single classical
$r$.
If $a\neq c$ then
$r({\mathfrak q}_{a,b,c})=(ae_{11}+be_{22}+ce_{33})\wedge e_{13}+
(a-c)e_{12}\wedge e_{23}$ is a classical $r$-matrix with carrier
${\mathfrak q}_{a,b,c}$.
We are unable to determine if this entire family consists
of boundary $r$-matrices, but there is a distinguished one-parameter
sub-family of ${\mathfrak q}_{a,b,c}$ consisting of boundary $r$-matrices.
If $x=e_{13}$ then a simple computation shows that
$\exp(t\ad x)\cdot r_{\lambda} = r_{\lambda}+t r({\mathfrak
q}_{-1-\lambda,
2\lambda, 1-\lambda})$ and thus $r({\mathfrak q}_{-1-\lambda,
2\lambda, 1-\lambda})$ is a boundary $r$-matrix for every
$\lambda$.
The quantization of any $r({\mathfrak q}_{a,b,c})$ is given by
$\exp(2t{\mathfrak q}_{a,b,c})$ which is
$$R({\mathfrak q}_{a,b,c})=\left(
\begin{matrix}
1 & 0 & a\,t & 0 & 0 & 0 &  - a\,t & 0 &  - a\,c\,t^{2} \\
0 & 1 & 0 & 0 & 0 & a\,t - c\,t & 0 &  - b\,t & 0 \\
0 & 0 & 1 & 0 & 0 & 0 & 0 & 0 &  - c\,t \\
0 & 0 & 0 & 1 & 0 & b\,t & 0 & c\,t - a\,t & 0 \\
0 & 0 & 0 & 0 & 1 & 0 & 0 & 0 & 0 \\
0 & 0 & 0 & 0 & 0 & 1 & 0 & 0 & 0 \\
0 & 0 & 0 & 0 & 0 & 0 & 1 & 0 & c\,t \\
0 & 0 & 0 & 0 & 0 & 0 & 0 & 1 & 0 \\
0 & 0 & 0 & 0 & 0 & 0 & 0 & 0 & 1
\end{matrix}\right).$$
For the cases $r({\mathfrak q}_{-1-\lambda,
2\lambda, 1-\lambda})$ Theorem \ref{thm3} also applies and the
resulting
boundary
$R$-matrix coincides with $R({\mathfrak q}_{-1-\lambda,
2\lambda, 1-\lambda})$.

As mentioned earlier, ${\mathfrak q}_{0,1,-1}$ and ${\mathfrak
q}_{1,1,-2}$
each carry an additional classical $r$. These are given by
$r'({\mathfrak q}_{0,1,-1})=(e_{22}-e_{33})\wedge e_{23} +
e_{12}\wedge e_{13}$
and $r'({\mathfrak
q}_{1,1,-2})=(e_{11}+e_{22}-2e_{33}+e_{23})\wedge e_{12}+
3e_{12}\wedge e_{23}$.
We do not know whether these are boundary $r$-matrices but,
fortunately, they
are easy to quantize. Their exponentials,
$\exp (2tr'({\mathfrak q}_{0,1,-1}))$ and
$\exp (2tr'({\mathfrak q}_{1,1,-2}))$, satisfy the quantum Yang-
Baxter equation.
These are given by
$$R'({\mathfrak q}_{0,1,-1})=
\left(
\begin{matrix}
1 & 0 & 0 & 0 & 0 & t & 0 &  - t & t^{2} \\
0 & 1 & 0 & 0 & 0 & 0 & 0 & 0 & 0 \\
0 & 0 & 1 & 0 & 0 & 0 & 0 & 0 & 0 \\
0 & 0 & 0 & 1 & 0 & 0 & 0 & 0 & 0 \\
0 & 0 & 0 & 0 & 1 & t & 0 &  - t & t^{2} \\
0 & 0 & 0 & 0 & 0 & 1 & 0 & 0 & t \\
0 & 0 & 0 & 0 & 0 & 0 & 1 & 0 & 0 \\
0 & 0 & 0 & 0 & 0 & 0 & 0 & 1 &  - t \\
0 & 0 & 0 & 0 & 0 & 0 & 0 & 0 & 1
\end{matrix}
\right)$$
and
$$R'({\mathfrak q}_{1,1,-2})=
\left(
\begin{matrix}
1 & 0 & t & 0 & 0 & 0 &  - t & 0 & 2\,t^{2} \\
0 & 1 & 0 & 0 & 0 & 3\,t & 0 &  - t &  - t \\
0 & 0 & 1 & 0 & 0 & 0 & 0 & 0 & 2\,t \\
0 & 0 & 0 & 1 & 0 & t & 0 &  - 3\,t & t \\
0 & 0 & 0 & 0 & 1 & 0 & 0 & 0 & 0 \\
0 & 0 & 0 & 0 & 0 & 1 & 0 & 0 & 0 \\
0 & 0 & 0 & 0 & 0 & 0 & 1 & 0 &  - 2\,t \\
0 & 0 & 0 & 0 & 0 & 0 & 0 & 1 & 0 \\
0 & 0 & 0 & 0 & 0 & 0 & 0 & 0 & 1
\end{matrix}
\right).$$
\item[(iv)] The two dimensional non-abelian Lie algebra has three
types of
embeddings in $\slt$. Each embedding carries a unique classical
$r$.
\begin{itemize}
\item[(a)]
 For any scalar $\lambda$ set
$${\mathfrak b}_{\lambda}=\left(
\begin{matrix}-1+\lambda & 0 & 0\\
0 & 1+\lambda & 0\\
0& 0 & -2\lambda
\end{matrix}\right) + \left(
\begin{matrix}0 & * & 0\\
0 & 0 & 0\\
0& 0 & 0
\end{matrix}\right).$$
Now $\exp(t\ad e_{12})\cdot r_{\lambda} =
r_{\lambda}+t ((-1+\lambda)e_{11}+(1+\lambda)e_{22}-2\lambda
e_{33})\wedge
e_{12}$
and so $r({\mathfrak b}_{\lambda})=
((-1+\lambda)e_{11}+(1+\lambda)e_{22}-2\lambda e_{33})\wedge
e_{12}$ is a
boundary
classical $r$-matrix.
Theorem \ref{thm3} then gives the boundary $R$-matrix
$$R({\mathfrak b}_{\lambda})=
\left(
\begin{matrix}
1 & \lambda\,t - t & 0 & t - \lambda\,t &  - \lambda^{2}\,t^{2} +
t^{2} & 0 & 0
& 0
 & 0 \\
0 & 1 & 0 & 0 &  - t - \lambda\,t & 0 & 0 & 0 & 0 \\
0 & 0 & 1 & 0 & 0 & 2\,\lambda\,t & 0 & 0 & 0 \\
0 & 0 & 0 & 1 & t + \lambda\,t & 0 & 0 & 0 & 0 \\
0 & 0 & 0 & 0 & 1 & 0 & 0 & 0 & 0 \\
0 & 0 & 0 & 0 & 0 & 1 & 0 & 0 & 0 \\
0 & 0 & 0 & 0 & 0 & 0 & 1 &  - 2\,\lambda\,t & 0 \\
0 & 0 & 0 & 0 & 0 & 0 & 0 & 1 & 0 \\
0 & 0 & 0 & 0 & 0 & 0 & 0 & 0 & 1
\end{matrix}
\right).$$

\item[(b)] Let
$${\mathfrak b}^{(0)}=*\left(
\begin{matrix}1&0&0\\0&0&0\\0&0&-1\end{matrix}\right)+
*\left(
\begin{matrix}0&1&0\\0&0&1\\0&0&0\end{matrix}\right).$$
Then $r({\mathfrak b}^{(0)})=  (e_{12}+e_{23})\wedge
(2e_{11}-2e_{33})$.
We are unable to determine if this is a boundary $r$-matrix.
Moreover,
$\exp(2t{\mathfrak b}^{(0)})$ does not satisfy the quantum Yang-Baxter
equation. Theorem \ref{thm4} however does apply here and
specializing the
universal $R$
gives
$$R({\mathfrak b}^{(0)})=
\left(
\begin{matrix}
1 &  - 2\,t & 2\,t^{2} & 2\,t & 0 & 0 & 2\,t^{2} & 0 & 0 \\
0 & 1 &  - 2\,t & 0 & 0 & 2\,t^{2} & 0 & 0 &  - 2\,t^{3} \\
0 & 0 & 1 & 0 & 0 &  - 2\,t & 0 & 0 & 2\,t^{2} \\
0 & 0 & 0 & 1 & 0 & 0 & 2\,t & 2\,t^{2} & 2\,t^{3} \\
0 & 0 & 0 & 0 & 1 & 0 & 0 & 0 & 0 \\
0 & 0 & 0 & 0 & 0 & 1 & 0 & 0 &  - 2\,t \\
0 & 0 & 0 & 0 & 0 & 0 & 1 & 2\,t & 2\,t^{2} \\
0 & 0 & 0 & 0 & 0 & 0 & 0 & 1 & 2\,t \\
0 & 0 & 0 & 0 & 0 & 0 & 0 & 0 & 1
\end{matrix}\right).$$
\item[(c)] Let
$${\mathfrak b}^{(1)}=*\left(
\begin{matrix}2&0&0\\0&-1&1\\0&0&-1\end{matrix}\right)+
\left(
\begin{matrix}0&0&*\\0&0&0\\0&0&0\end{matrix}\right).$$
The classical $r$ here is
$r({\mathfrak b}^{(1)})=  e_{13}\wedge (2e_{11}-e_{22}-
e_{33}+e_{23})$
but just like case (b), we are unable to determine if this is a
boundary
$r$-matrix but specializing the universal $R$ of Theorem
\ref{thm4}
for this representation gives
the $R$-matrix
$$R({\mathfrak b}^{(1)})=
\left(
\begin{matrix}
1 & 0 &  - 2\,t & 0 & 0 & 0 & 2\,t & 0 & 2\,t^{2} \\
0 & 1 & 0 & 0 & 0 & 0 & 0 &  - t & t \\
0 & 0 & 1 & 0 & 0 & 0 & 0 & 0 &  - t \\
0 & 0 & 0 & 1 & 0 & t & 0 & 0 &  - t \\
0 & 0 & 0 & 0 & 1 & 0 & 0 & 0 & 0 \\
0 & 0 & 0 & 0 & 0 & 1 & 0 & 0 & 0 \\
0 & 0 & 0 & 0 & 0 & 0 & 1 & 0 & t \\
0 & 0 & 0 & 0 & 0 & 0 & 0 & 1 & 0 \\
0 & 0 & 0 & 0 & 0 & 0 & 0 & 0 & 1
\end{matrix}\right).$$

\begin{rem} In \cite{CP} and \cite{S} it is claimed that there
are four embeddings of the two-dimensional non-abelian Lie
algebra in $\slt$. However, the algebra of type ${\mathfrak
b}^{(2)}$, in \cite{CP} is three-dimensional and the algebra of
type $C_{1/3}^{1,1}$ in \cite{S} is conjugate to type ${\mathfrak
b}^{(c)}$ discussed in item (a).
\end{rem}

\item[(v)] There are three embeddings of the two-dimensional
abelian Lie algebra in $\slt$.

\begin{itemize}
\item[(a)] The Cartan subalgebra
$${\mathfrak h}=\left(
\begin{matrix}*&0&0\\0&*&0\\0&0&*\end{matrix}\right)$$
carries the unique classical $r$-matrix
$r({\mathfrak h})=(e_{11}-e_{22})\wedge (e_{22}-e_{33}).$
It is a boundary $r$-matrix since
$\lim _{\lambda \to \infty}(1/\lambda)r_{\lambda}= r({\mathfrak
h}).$
The modified $R$-matrix which quantizes
$r_{\lambda}$ is
$\exp (t\lambda r({\mathfrak h}))\cdot \exp(2t\gamma)\cdot
\exp (t\lambda r({\mathfrak h})).$ Denote this modified $R$-
matrix by $Q$.
Then if $t=m/\lambda$, we
have that $\lim _{\lambda \to \infty}Q=\exp (2m r({\mathfrak
h}))$ which,
being a diagonal matrix, satisfies the quantum Yang-Baxter
equation.
Hence
$$R({\mathfrak h})=
\left(
\begin{matrix}
1 & 0 &  0 & 0 & 0 & 0 & 0 & 0 & 0 \\
0 & e^{m} & 0 & 0 & 0 & 0 & 0 &  0 &  \\
0 & 0 & e^{-m} & 0 & 0 & 0 & 0 & 0 &  0 \\
0 & 0 & 0 & e^{-m} & 0 & 0 & 0 & 0 &  0 \\
0 & 0 & 0 & 0 & 1 & 0 & 0 & 0 & 0 \\
0 & 0 & 0 & 0 & 0 & e^{m} & 0 & 0 & 0 \\
0 & 0 & 0 & 0 & 0 & 0 & e^{m} & 0 & 0 \\
0 & 0 & 0 & 0 & 0 & 0 & 0 & e^{-m} & 0 \\
0 & 0 & 0 & 0 & 0 & 0 & 0 & 0 & 1
\end{matrix}\right)$$
is a boundary $R$-matrix.
\item[(b)] The abelian subalgebra
$${\mathfrak h}^{(1)}=*\left(
\begin{matrix}1&0&0\\0&1&0\\0&0&-2\end{matrix}\right)+
*\left(
\begin{matrix}0&1&0\\0&0&0\\0&0&0\end{matrix}\right)$$
is the carrier of $r({\mathfrak
h}^{(1)})=(e_{11}+e_{22}-2e_{33})\wedge
e_{12}.$
It is a boundary $r$-matrix because
if $r({\mathfrak b}_{\lambda})$ is the boundary $r$-matrix
associated to the
non-abelian Lie algebra ${\mathfrak b}_{\lambda}$,
then $\lim _{\lambda \to \infty} (1/\lambda)r({\mathfrak
b}_{\lambda})=
r({\mathfrak h}^{(1)})$ and so it is a limiting case of boundary
classical
$r$-matrices and so must be one itself.
The same argument works on the quantum level, giving the boundary
$R$-matrix
$$R({\mathfrak h}^{(1)})=
\left(
\begin{matrix}
1 & t & 0 &  - t &  - t^{2} & 0 & 0 & 0 & 0 \\
0 & 1 & 0 & 0 &  - t & 0 & 0 & 0 & 0 \\
0 & 0 & 1 & 0 & 0 & 2\,t & 0 & 0 & 0 \\
0 & 0 & 0 & 1 & t & 0 & 0 & 0 & 0 \\
0 & 0 & 0 & 0 & 1 & 0 & 0 & 0 & 0 \\
0 & 0 & 0 & 0 & 0 & 1 & 0 & 0 & 0 \\
0 & 0 & 0 & 0 & 0 & 0 & 1 &  - 2\,t & 0 \\
0 & 0 & 0 & 0 & 0 & 0 & 0 & 1 & 0 \\
0 & 0 & 0 & 0 & 0 & 0 & 0 & 0 & 1
\end{matrix}\right).$$

\item[(c)] The only other abelian subalgebras of $\slt$ are
$${\mathfrak h}^{(\lambda,1)}=*\left(
\begin{matrix}0&1&0\\0&0&\lambda\\0&0&0\end{matrix}\right)
+\left(
\begin{matrix}0&0&*\\0&0&0\\0&0&0\end{matrix}\right).
$$
Up to conjugacy however, there are only two algebras in this
family,
${\mathfrak h}^{(1,1)}$ and ${\mathfrak h}^{(0,1)}$. In any case,
the associated classical $r$-matrix is
$r({\mathfrak h}^{(\lambda,1)})=(e_{12}+\lambda e_{23})\wedge
e_{13}$.
This is also a limiting case of boundary solutions. Recall that
the Lie algebra ${\mathfrak q}_{-1,0,1}$ is in the distinguished
sub-family of
${\mathfrak q}_{a,b,c}$ which carries boundary solutions
 Set $x=e_{12}+ \lambda e_{13}$. Then
$\exp (t\ad x)\cdot r({\mathfrak q}_{-1,0,1})= r({\mathfrak
q}_{-1,
0,1 })-2t r({\mathfrak h}^{(\lambda,1)})$
and so $r({\mathfrak h}^{(\lambda,1)})$ is a boundary $r$-matrix.
The associated boundary $R$-matrix is found in the same way,
specifically
it is
$$R({\mathfrak h}^{(\lambda,1)})=
\left(
\begin{matrix}
1 & 0 & 0 & 0 & 0 & t & 0 &  - t & 0 \\
0 & 1 & 0 & 0 & 0 & 0 & 0 & 0 &  - \lambda \,t \\
0 & 0 & 1 & 0 & 0 & 0 & 0 & 0 & 0 \\
0 & 0 & 0 & 1 & 0 & 0 & 0 & 0 & \lambda \,t \\
0 & 0 & 0 & 0 & 1 & 0 & 0 & 0 & 0 \\
0 & 0 & 0 & 0 & 0 & 1 & 0 & 0 & 0 \\
0 & 0 & 0 & 0 & 0 & 0 & 1 & 0 & 0 \\
0 & 0 & 0 & 0 & 0 & 0 & 0 & 1 & 0 \\
0 & 0 & 0 & 0 & 0 & 0 & 0 & 0 & 1
\end{matrix}\right).$$
\end{itemize}
\end{itemize}

     \end{enumerate}

\end{document}